\documentclass[12pt,reqno]{amsart}
\advance\textheight by \topskip
\allowdisplaybreaks
\title[Elementary proof of a series evaluation]{An elementary proof of a series evaluation in terms of harmonic numbers}

\author[H.~Prodinger]{Helmut Prodinger}
\address{Helmut Prodinger\\
Department of Mathematics\\
University of Stellenbosch\\
7602 Stellenbosch\\
South Africa}
\email{hproding@sun.ac.za}

\date{May 20, 2008}

\begin{document}

\maketitle

For positive integers $j$, consider
\begin{equation*}
S(j)=\sum_{n\ge1}\frac1{2^{2n-1}(2n-1)}\sum_k \binom{2n-1}{k}\frac1{k-j-n+\frac12}.
\end{equation*}

This quantity arose in \cite{LySt03} and was subsequently evaluated in \cite{LyPaRi02}.
Further proofs of the final formula (\ref{answer}) were given in \cite{Krattenthaler04,ChDo05}.

Here, we give an extremely short and simple proof.
\medskip

For our analysis, it is better to replace the inner summation index $k$ by $2n-1-k$, and consider
\begin{equation*}
T(j)=\sum_{n\ge1}\frac1{2^{2n-1}(2n-1)}\sum_k \binom{2n-1}{k}\frac1{k+j-n+\frac12},
\end{equation*}
then $S(j)=-T(j)$.

We start from the obvious formula
\begin{equation*}
\int_0^1t^{x-1}(1+t)^{2n-1}dt=\sum_k \binom{2n-1}{k}\frac1{k+x},
\end{equation*}
substitute $x=j-n+\frac12$ and sum:
\begin{align*}
T(j)&=\int_0^1t^{j-\frac12}(1+t)^{-1}\sum_{n\ge1}\frac1{2^{2n-1}(2n-1)}t^{-n}(1+t)^{2n}dt\\
&=\int_0^1t^{j-\frac12}(1+t)^{-1}\frac{1+t}{\sqrt t}\log\frac{1+\sqrt t}{1-\sqrt t}dt\\
&=2\int_0^1w^{2j-1}\log\frac{1+w}{1-w}dw\\
&=4\int_0^1w^{2j-1}\sum_{k\ge1}\frac{w^{2k-1}}{2k-1}dw\\
&=4\sum_{k\ge1}\frac{1}{(2k-1)(2k+2j-1)}\\
&=\frac2j\sum_{k\ge1}\Big(\frac{1}{2k-1}-\frac{1}{2k+2j-1}\Big)\\
&=\frac2j\sum_{k=1}^{j}\frac{1}{2k-1}.
\end{align*}
Hence
\begin{equation}\label{answer}
S(j)=-\frac2j\sum_{k=1}^{j}\frac{1}{2k-1}.
\end{equation}
It is not hard to show that it is allowed to interchange integration and summation.

\bibliographystyle{plain}

\end{document}